# Synchronal and cyclic algorithms for fixed point problems and variational inequality problems in Banach spaces

Abba Auwalu[1*], Lawan Bulama Mohammed[2] and Afis Saliu[3]

*Correspondence:
abbaauwalu@yahoo.com
[1]Department of Mathematical Sciences, College of Remedial and Advanced Studies, P.M.B. 048, Kafin Hausa, Jigawa, Nigeria
Full list of author information is available at the end of the article

**Abstract**

In this paper, we study synchronal and cyclic algorithms for finding a common fixed point $x^*$ of a finite family of strictly pseudocontractive mappings, which solve the variational inequality

$$\langle (\gamma f - \mu G) x^*, j_q(x - x^*) \rangle \leq 0, \quad \forall x \in \bigcap_{i=1}^{N} F(T_i),$$

where $f$ is a contraction mapping, $G$ is an $\eta$-strongly accretive and $L$-Lipschitzian operator, $N \geq 1$ is a positive integer, $\gamma, \mu > 0$ are arbitrary fixed constants, and $\{T_i\}_{i=1}^{N}$ are $N$-strict pseudocontractions. Furthermore, we prove strong convergence theorems of such iterative algorithms in a real $q$-uniformly smooth Banach space. The results presented extend, generalize and improve the corresponding results recently announced by many authors.

## 1 Introduction

Let $E$ be a real Banach space, and let $E^*$ be the dual of $E$. For some real number $q$ ($1 < q < \infty$), the *generalized duality mapping* $J_q : E \to 2^{E^*}$ is defined by

$$J_q(x) = \{x^* \in E^* : \langle x, x^* \rangle = \|x\|^q, \|x^*\| = \|x\|^{q-1}\}, \quad \forall x \in E, \tag{1.1}$$

where $\langle \cdot, \cdot \rangle$ denotes the duality pairing between elements of $E$ and those of $E^*$. In particular, $J = J_2$ is called the *normalized duality mapping* and $J_q(x) = \|x\|^{q-2} J_2(x)$ for $x \neq 0$. If $E$ is a real Hilbert space, then $J = I$, where $I$ is the identity mapping. It is well known that if $E$ is smooth, then $J_q$ is single-valued, which is denoted by $j_q$.

Let $C$ be a nonempty closed convex subset of $E$, and let $G : E \to E$ be a nonlinear map. Then the *variational inequality problem* with respect to $C$ and $G$ is to find a point $x^* \in C$ such that

$$\langle Gx^*, j_q(x - x^*) \rangle \geq 0, \quad \forall x \in C \quad \text{and} \quad j_q(x - x^*) \in J_q(x - x^*). \tag{1.2}$$

We denote by $VI(G, C)$ the set of solutions of this variational inequality problem.

If $E = H$, a real Hilbert space, the variational inequality problem reduces to the following: Find a point $x^* \in C$ such that

$$\langle Gx^*, x - x^* \rangle \geq 0, \quad \forall x \in C. \tag{1.3}$$

A mapping $T : E \to E$ is said to be a *contraction* if, for some $\alpha \in [0, 1)$,

$$\|Tx - Ty\| \leq \alpha \|x - y\|, \quad \forall x, y \in E. \tag{1.4}$$

The map $T$ is said to be *nonexpansive* if

$$\|Tx - Ty\| \leq \|x - y\|, \quad \forall x, y \in E. \tag{1.5}$$

The map $T$ is said to be *L-Lipschitzian* if there exists $L > 0$ such that

$$\|Tx - Ty\| \leq L \|x - y\|, \quad \forall x, y \in E. \tag{1.6}$$

A point $x \in E$ is called a fixed point of the map $T$ if $Tx = x$. We denote by $F(T)$ the set of all fixed points of the mapping $T$, that is,

$$F(T) = \{x \in C : Tx = x\}.$$

We assume that $F(T) \neq \emptyset$ in the sequel. It is well known that $F(T)$ above is closed and convex (see, *e.g.*, Goebel and Kirk [1]).

An operator $F : E \to E$ is said to be accretive if $\forall x, y \in E$, there exists $j_q(x - y) \in J_q(x - y)$ such that

$$\langle Fx - Fy, j_q(x - y) \rangle \geq 0. \tag{1.7}$$

For some positive real numbers $\eta$, $\lambda$, the mapping $F$ is said to be $\eta$-*strongly accretive* if for any $x, y \in E$, there exists $j_q(x - y) \in J_q(x - y)$ such that

$$\langle Fx - Fy, j_q(x - y) \rangle \geq \eta \|x - y\|^q, \tag{1.8}$$

and it is called $\lambda$-*strictly pseudocontractive* if

$$\langle Fx - Fy, j_q(x - y) \rangle \leq \|x - y\|^q - \lambda \|x - y - (Fx - Fy)\|^q. \tag{1.9}$$

It is clear that (1.9) is equivalent to the following:

$$\langle (I - F)x - (I - F)y, j_q(x - y) \rangle \geq \lambda \|x - y - (Fx - Fy)\|^q, \tag{1.10}$$

where $I$ denotes the identity operator.

In Hilbert spaces, *accretive* operators are called *monotone* where inequality (1.7) holds with $j_q$ replaced by the identity map of $H$.

A bounded linear operator A on H is called strongly positive with coefficient $\gamma$ if there is a constant $\gamma > 0$ with the property

$$\langle Ax, x \rangle \geq \gamma \|x\|^2, \quad \forall x \in H.$$

Let K be a nonempty closed convex and bounded subset of a Banach space E, and let the diameter of K be defined by $d(K) := \sup\{\|x - y\| : x, y \in K\}$. For each $x \in K$, let $r(x, K) := \sup\{\|x - y\| : y \in K\}$, and let $r(K) := \inf\{r(x, K) : x \in K\}$ denote the Chebyshev radius of $K$ relative to itself. The normal structure coefficient $N(E)$ of $E$ (see, *e.g.*, [2]) is defined by $N(E) := \inf\{\frac{d(K)}{r(K)} : d(K) > 0\}$. A space $E$, such that $N(E) > 1$, is said to have a uniform normal structure.

It is known that all uniformly convex and uniformly smooth Banach spaces have a uniform normal structure (see, *e.g.*, [3, 4]).

Let $\mu$ be a continuous linear functional on $l^\infty$ and $(a_0, a_1, \ldots) \in l^\infty$. We write $\mu_n(a_n)$ instead of $\mu((a_0, a_1, \ldots))$. We call $\mu$ a Banach limit if $\mu$ satisfies $\|\mu\| = \mu_n(1) = 1$ and $\mu_n(a_{n+1}) = \mu_n(a_n)$ for all $(a_0, a_1, \ldots) \in l^\infty$. If $\mu$ is a Banach limit, then

$$\liminf_{n \to \infty} a_n \leq \mu_n(a_n) \leq \limsup_{n \to \infty} a_n$$

for all $\{a_n\} \in l^\infty$ (see, *e.g.*, [3, 5]).

Let $S = \{x \in E : \|x\| = 1\}$ denote the unit sphere of a real Banach space $E$.

The space $E$ is said to have a Gâteaux differentiable norm if the limit

$$\lim_{t \to 0} \frac{\|x + ty\| - \|x\|}{t} \tag{1.11}$$

exists for each $x, y \in S$. In this case, $E$ is called smooth. $E$ is said to be uniformly smooth if the limit (1.11) exists and is attained uniformly in $x, y \in S$. $E$ is said to have a uniformly Gâteaux differentiable norm if, for any $y \in S$, the limit (1.11) exists uniformly for all $x \in S$.

The modulus of smoothness of $E$, with $\dim E \geq 2$, is a function $\rho_E : [0, \infty) \to [0, \infty)$ defined by

$$\rho_E(\tau) = \sup\left\{\frac{\|x + y\| + \|x - y\|}{2} - 1 : \|x\| = 1, \|y\| \leq \tau\right\}.$$

A Banach space $E$ is said to be uniformly smooth if $\lim_{t \to 0^+} \frac{\rho_E(t)}{t} = 0$, and for $q > 1$, $E$ is said to be $q$-uniformly smooth if there exists a fixed constant $c > 0$ such that $\rho_E(t) \leq ct^q$, $t > 0$.

It is well known (see, *e.g.*, [6]) that Hilbert spaces, $L_p$ (or $l_p$) spaces ($1 < p < \infty$) and Sobolev spaces, $W_m^p$ ($1 < p < \infty$) are all uniformly smooth. More precisely, Hilbert spaces are 2-uniformly smooth, while

$$L_p \text{ (or } l_p\text{) or } W_m^p \text{ spaces are } \begin{cases} 2\text{-uniformly smooth if } 2 \leq p < \infty, \\ p\text{-uniformly smooth if } 1 < p \leq 2. \end{cases}$$

Also, it is well known (see, *e.g.*, [7]) that $q$-uniformly smooth Banach spaces have a uniformly Gâteaux differentiable norm.

The variational inequality problem was initially introduced and studied by Stampacchia [8] in 1964. In the recent years, variational inequality problems have been extended to study a large variety of problems arising in structural analysis, economics and optimization. Thus, the problem of solving a variational inequality of the form (1.2) has been intensively studied by numerous authors. Iterative methods for approximating fixed points of nonexpansive mappings and their generalizations, which solve some variational inequality problems, have been studied by a number of authors (see, for example, [9–13] and the references therein).

Let $H$ be a real Hilbert space. In 2001, Yamada [13] proposed a hybrid steepest descent method for solving variational inequality as follows: Let $x_0 \in H$ be chosen arbitrarily and define a sequence $\{x_n\}$ by

$$x_{n+1} = Tx_n - \mu\lambda_n F(Tx_n), \quad n \geq 0, \tag{1.12}$$

where $T$ is a nonexpansive mapping on $H$, $F$ is $L$-Lipschitzian and $\eta$-strongly monotone with $L > 0$, $\eta > 0$, $0 < \mu < 2\eta/L^2$. If $\{\lambda_n\}$ is a sequence in $(0,1)$ satisfying the following conditions:

(C1) $\lim_{n\to\infty} \lambda_n = 0$,
(C2) $\sum_{n=0}^{\infty} \lambda_n = \infty$,
(C3) either $\sum_{n=1}^{\infty} |\lambda_{n+1} - \lambda_n| < \infty$ or $\lim_{n\to\infty} \frac{\lambda_{n+1}}{\lambda_n} = 1$,

then he proved that the sequence $\{x_n\}$ converges strongly to the unique solution of the variational inequality

$$\langle F\tilde{x}, x - \tilde{x}\rangle \geq 0, \quad \forall x \in F(T).$$

Besides, he also proposed the cyclic algorithm

$$x_{n+1} = T^{\lambda_n} x_n = (I - \mu\lambda_n F)T_{[n]}x_n,$$

where $T_{[n]} = T_{n(\mathrm{mod}\, N)}$; he also proved strong convergence theorems for the cyclic algorithm.

In 2006, Marino and Xu [10] considered the following general iterative method: Starting with an arbitrary initial point $x_0 \in H$, define a sequence $\{x_n\}$ by

$$x_{n+1} = \alpha_n \gamma f(x_n) + (I - \alpha_n A)Tx_n, \quad n \geq 0, \tag{1.13}$$

where $T$ is a nonexpansive mapping of $H$, $f$ is a contraction, $A$ is a linear bounded strongly positive operator, and $\{\alpha_n\}$ is a sequence in $(0,1)$ satisfying the following conditions:

(M1) $\lim_{n\to\infty} \alpha_n = 0$;
(M2) $\sum_{n=0}^{\infty} \alpha_n = \infty$;
(M3) $\sum_{n=1}^{\infty} |\alpha_{n+1} - \alpha_n| < \infty$ or $\lim_{n\to\infty} \frac{\alpha_{n+1}}{\alpha_n} = 1$.

They proved that the sequence $\{x_n\}$ converges strongly to a fixed point $\tilde{x}$ of $T$, which solves the variational inequality

$$\langle(\gamma f - A)\tilde{x}, x - \tilde{x}\rangle \leq 0, \quad \forall x \in F(T).$$

In 2010, Tian [11] combined the iterative method (1.13) with Yamada's iterative method (1.12) and considered the following general iterative method:

$$x_{n+1} = \alpha_n \gamma f(x_n) + (I - \mu \alpha_n F) T x_n, \quad n \geq 0, \tag{1.14}$$

where $T$ is a nonexpansive mapping on $H$, $f$ is a contraction, $F$ is $k$-Lipschitzian and $\eta$-strongly monotone with $k > 0$, $\eta > 0$, $0 < \mu < 2\eta/k^2$. He proved that if the sequence $\{\alpha_n\}$ of parameters satisfies conditions (M1)-(M3), then the sequence $\{x_n\}$ generated by (1.14) converges strongly to a fixed point $\tilde{x}$ of $T$, which solves the variational inequality

$$\langle (\gamma f - \mu F)\tilde{x}, x - \tilde{x} \rangle \leq 0, \quad \forall x \in F(T).$$

Very recently, in 2011, Tian and Di [12] studied two algorithms, based on Tian's [11] general iterative algorithm, and proved the following theorems.

**Theorem 1.1** (Synchronal algorithm) *Let $H$ be a real Hilbert space, and let $T_i : H \to H$ be a $k_i$-strictly pseudocontraction for some $k_i \in (0, 1)$ ($i = 1, 2, \ldots, N$) such that $\bigcap_{i=1}^{N} F(T_i) \neq \emptyset$, and $f$ be a contraction with coefficient $\beta \in (0, 1)$ and $\lambda_i$ be a positive constant such that $\sum_{i=1}^{N} \lambda_i = 1$. Let $G : H \to H$ be an $\eta$-strongly monotone and L-Lipschitzian operator with $L > 0$, $\eta > 0$. Assume that $0 < \mu < 2\eta/L^2$, $0 < \gamma < \mu(\eta - \frac{\mu L^2}{2})/\beta = \tau/\beta$. Let $x_0 \in H$ be chosen arbitrarily, and let $\{\alpha_n\}$ and $\{\beta_n\}$ be sequences in $(0, 1)$ satisfying the following conditions*:

(N1) $\lim_{n \to \infty} \alpha_n = 0$, $\sum_{n=0}^{\infty} \alpha_n = \infty$;
(N2) $\sum_{n=1}^{\infty} |\alpha_{n+1} - \alpha_n| < \infty$, $\sum_{n=1}^{\infty} |\beta_{n+1} - \beta_n| < \infty$;
(N3) $0 < \max k_i \leq \beta_n < a < 1$, $\forall n \geq 0$.

*Let $\{x_n\}$ be a sequence defined by the composite process*

$$\begin{cases} T^{\beta_n} = \beta_n I + (1 - \beta_n) \sum_{i=1}^{N} \lambda_i T_i, \\ x_{n+1} = \alpha_n \gamma f(x_n) + (I - \alpha_n \mu G) T^{\beta_n} x_n, \quad n \geq 0. \end{cases}$$

*Then $\{x_n\}$ converges strongly to a common fixed point of $\{T_i\}_{i=1}^{N}$, which solves the variational inequality*

$$\langle (\gamma f - \mu G) x^*, x - x^* \rangle \leq 0, \quad \forall x \in \bigcap_{i=1}^{N} F(T_i). \tag{1.15}$$

**Theorem 1.2** (Cyclic algorithm) *Let $H$ be a real Hilbert space, and let $T_i : H \to H$ be a $k_i$-strictly pseudocontraction for some $k_i \in (0, 1)$ ($i = 1, 2, \ldots, N$) such that $\bigcap_{i=1}^{N} F(T_i) \neq \emptyset$, and $f$ be a contraction with coefficient $\beta \in (0, 1)$. Let $G : H \to H$ be an $\eta$-strongly monotone and L-Lipschitzian operator with $L > 0$, $\eta > 0$. Assume that $0 < \mu < 2\eta/L^2$, $0 < \gamma < \mu(\eta - \frac{\mu L^2}{2})/\beta = \tau/\beta$. Let $x_0 \in H$ be chosen arbitrarily, and let $\{\alpha_n\}$ and $\{\beta_n\}$ be sequences in $(0, 1)$ satisfying the following conditions*:

(N1′) $\lim_{n \to \infty} \alpha_n = 0$, $\sum_{n=0}^{\infty} \alpha_n = \infty$;
(N2′) $\sum_{n=1}^{\infty} |\alpha_{n+1} - \alpha_n| < \infty$ or $\lim_{n \to \infty} \frac{\alpha_n}{\alpha_{n+N}} = 1$;
(N3′) $\beta_{[n]} \in [k, 1)$, $\forall n \geq 0$, where $k = \max\{k_i : 1 \leq i \leq N\}$.

*Let $\{x_n\}$ be a sequence defined by the composite process*

$$\begin{cases} A_{[n]} = \beta_{[n]}I + (1 - \beta_{[n]})T_{[n]}, \\ x_{n+1} = \alpha_n \gamma f(x_n) + (I - \alpha_n \mu G)A_{[n+1]}x_n, \quad n \geq 0, \end{cases}$$

*where $T_{[n]} = T_i$, with $i = n(\mod N)$, $1 \leq i \leq N$, namely $T_{[n]}$ is one of $T_1, T_2, \ldots, T_N$ cyclically. Then $\{x_n\}$ converges strongly to a common fixed point of $\{T_i\}_{i=1}^N$, which solves the variational inequality* (1.15).

In this paper, we study the synchronal and cyclic algorithms for finding a common fixed point $x^*$ of finite strictly pseudocontractive mappings, which solves the variational inequality

$$\langle (\gamma f - \mu G)x^*, j_q(x - x^*) \rangle \leq 0, \quad \forall x \in \bigcap_{i=1}^N F(T_i), \tag{1.16}$$

where $f$ is a contraction mapping, $G$ is an $\eta$-strongly accretive and $L$-Lipschitzian operator, $N \geq 1$ is a positive integer, $\gamma, \mu > 0$ are arbitrary fixed constants, and $\{T_i\}_{i=1}^N$ are $N$-strict pseudocontractions defined on a closed convex subset $C$ of a real $q$-uniformly smooth Banach space $E$ whose norm is uniformly Gâteaux differentiable.

Let $T$ be defined by

$$T := \sum_{i=1}^N \lambda_i T_i,$$

where $\lambda_i > 0$ such that $\sum_{i=1}^N \lambda_i = 1$. We will show that a sequence $\{x_n\}$ generated by the following *synchronal algorithm*:

$$\begin{cases} x_0 = x \in C \quad \text{chosen arbitrarily}, \\ T^{\beta_n} = \beta_n I + (1 - \beta_n)\sum_{i=1}^N \lambda_i T_i, \\ x_{n+1} = \alpha_n \gamma f(x_n) + (I - \alpha_n \mu G)T^{\beta_n} x_n, \quad n \geq 0, \end{cases} \tag{1.17}$$

converges strongly to a solution of problem (1.16).

Another approach to problem (1.16) is the *cyclic algorithm*. For each $i = 1, \ldots, N$, let $A_i = \beta_i I + (1 - \beta_i)T_i$, where the constant $\beta_i$ satisfies $0 < k_i < \beta_i < 1$. Beginning with $x_0 \in C$, define a sequence $\{x_n\}$ cyclically by

$$x_1 = \alpha_0 \gamma f(x_0) + (I - \alpha_0 \mu G)(A_1 x_0),$$

$$x_2 = \alpha_1 \gamma f(x_1) + (I - \alpha_1 \mu G)(A_2 x_1),$$

$$\vdots$$

$$x_N = \alpha_{N-1} \gamma f(x_{N-1}) + (I - \alpha_{N-1} \mu G)(A_N x_{N-1}),$$

$$x_{N+1} = \alpha_N \gamma f(x_N) + (I - \alpha_N \mu G)(A_1 x_N),$$

$$\vdots$$

Indeed, the algorithm can be written in a compact form as follows:

$$\begin{cases} x_0 = x \in C \quad \text{chosen arbitrarily,} \\ A_{[n]} = \beta_{[n]} I + (1 - \beta_{[n]}) T_{[n]}, \\ x_{n+1} = \alpha_n \gamma f(x_n) + (I - \alpha_n \mu G) A_{[n+1]} x_n, \quad n \geq 0, \end{cases} \quad (1.18)$$

where $T_{[n]} = T_i$, with $i = n \pmod{N}$, $1 \leq i \leq N$, namely $T_{[n]}$ is one of $T_1, T_2, \ldots, T_N$ cyclically. We will show that (1.18) is also strongly convergent to a solution of problem (1.16) if the sequences $\{\alpha_n\}$ and $\{\beta_n\}$ of parameters are appropriately chosen.

Motivated by the results of Tian and Di [12], in this paper we aim to continue the study of fixed point problems and prove new theorems for the solution of variational inequality problems in the framework of a real Banach space, which is much more general than that of Hilbert.

Throughout this research work, we will use the following notations:

1. $\rightharpoonup$ for weak convergence and $\to$ for strong convergence.
2. $\omega_\omega(x_n) = \{x : \exists x_{n_j} \rightharpoonup x\}$ denotes the weak $\omega$-limit set of $\{x_n\}$.

## 2 Preliminaries

In the sequel we shall make use of the following lemmas.

**Lemma 2.1** (Zhang and Guo [14]) *Let $C$ be a nonempty closed convex subset of a real Banach space $E$. Given an integer $N \geq 1$, for each $1 \leq i \leq N$, $T_i : C \to C$ is a $\lambda_i$-strict pseudo-contraction for some $\lambda_i \in [0,1)$ such that $\bigcap_{i=1}^{N} F(T_i) \neq \emptyset$. Assume that $\{\gamma_i\}_{i=1}^{N}$ is a sequence of positive numbers such that $\sum_{i=1}^{N} \gamma_i = 1$. Then $\sum_{i=1}^{N} \gamma_i T_i$ is a $\lambda$-strict pseudocontraction with $\lambda := \min\{\lambda_i : 1 \leq i \leq N\}$, and*

$$F\left(\sum_{i=1}^{N} \gamma_i T_i\right) = \bigcap_{i=1}^{N} F(T_i).$$

**Lemma 2.2** (Zhou [15]) *Let $E$ be a uniformly smooth real Banach space, and let $C$ be a nonempty closed convex subset of $E$. Let $T : C \to C$ be a $k$-strict pseudocontraction. Then $(I - T)$ is demiclosed at zero. That is, if $\{x_n\} \subset C$ satisfies $x_n \rightharpoonup x$ and $x_n - Tx_n \to 0$, as $n \to \infty$, then $Tx = x$.*

**Lemma 2.3** (Petryshyn [16]) *Let $E$ be a real Banach space, and let $J_q : E \to 2^{E^*}$ be the generalized duality mapping. Then, for any $x, y \in E$ and $j_q(x+y) \in J_q(x+y)$,*

$$\|x + y\|^q \leq \|x\|^q + q\langle y, j_q(x+y) \rangle.$$

**Lemma 2.4** (Lim and Xu [4]) *Suppose that $E$ is a Banach space with a uniform normal structure, $K$ is a nonempty bounded subset of $E$, and let $T : K \to K$ be a uniformly $k$-Lipschitzian mapping with $k < N(E)^{\frac{1}{2}}$. Suppose also that there exists a nonempty bounded closed convex subset $C$ of $K$ with the following property* (P):

$$x \in C \quad \text{implies} \quad \omega_\omega(x) \subset C, \quad (P)$$

where $\omega_\omega(x)$ is the $\omega$-limit set of $T$ at $x$, i.e., the set

$$\left\{y \in E : y = \text{weak} - \lim_j T^{n_j}x \text{ for some } n_j \to \infty\right\}.$$

*Then $T$ has a fixed point in $C$.*

**Lemma 2.5** (Xu [17]) *Let $q > 1$, and let $E$ be a real $q$-uniformly smooth Banach space, then there exists a constant $d_q > 0$ such that for all $x, y \in E$ and $j_q(x) \in J_q(x)$,*

$$\|x + y\|^q \leq \|x\|^q + q\langle y, j_q(x)\rangle + d_q\|y\|^q.$$

**Lemma 2.6** *Let $E$ be a real $q$-uniformly smooth Banach space with constant $d_q > 0$, $q > 1$, and let $C$ be a nonempty closed convex subset of $E$. Let $F : C \to C$ be an $\eta$-strongly accretive and $L$-Lipschitzian operator with $L > 0$, $\eta > 0$. Assume that $0 < \mu < (\frac{q\eta}{d_q L^q})^{\frac{1}{q-1}}$, $\tau = \mu(\eta - \frac{d_q\mu^{q-1}L^q}{q})$ and $t \in (0, \min\{1, \frac{1}{\tau}\})$. Then, for any $x, y \in C$, the following inequality holds:*

$$\|(I - \mu tF)x - (I - \mu tF)y\| \leq (1 - t\tau)\|x - y\|.$$

*That is, $(I - \mu tF)$ is a contraction with coefficient $(1 - t\tau)$.*

*Proof* For any $x, y \in C$, we have, by Lemma 2.5, (1.6) and (1.8),

$$\begin{aligned}
\|(I - \mu tF)x - (I - \mu tF)y\|^q &= \|(x - y) - \mu t(Fx - Fy)\|^q \\
&\leq \|x - y\|^q - q\mu t\langle Fx - Fy, j_q(x - y)\rangle + d_q\mu^q t^q\|Fx - Fy\|^q \\
&\leq \|x - y\|^q - q\mu t\eta\|x - y\|^q + d_q\mu^q t^q L^q\|x - y\|^q \\
&\leq \left[1 - t\mu(q\eta - d_q\mu^{q-1}L^q)\right]\|x - y\|^q \\
&= \left[1 - qt\mu\left(\eta - \frac{d_q\mu^{q-1}L^q}{q}\right)\right]\|x - y\|^q \\
&\leq \left[1 - t\mu\left(\eta - \frac{d_q\mu^{q-1}L^q}{q}\right)\right]^q\|x - y\|^q \\
&= (1 - t\tau)^q\|x - y\|^q.
\end{aligned}$$

From $0 < \mu < (\frac{q\eta}{d_q L^q})^{\frac{1}{q-1}}$, $q > 1$ and $t \in (0, \min\{1, \frac{1}{\tau}\})$, we have $(1 - t\tau) \in (0, 1)$. It then follows that

$$\|(I - \mu tF)x - (I - \mu tF)y\| \leq (1 - t\tau)\|x - y\|. \qquad \square$$

**Lemma 2.7** *Let $E$ be a real $q$-uniformly smooth Banach space with constant $d_q$, $q > 1$, and let $C$ be a nonempty closed convex subset of $E$. Suppose that $T : C \to C$ is a $\lambda$-strict pseudocontraction such that $F(T) \neq \emptyset$. For any $\alpha \in (0, 1)$, we define $T_\alpha : C \to E$ by $T_\alpha x = \alpha x + (1 - \alpha)Tx$ for each $x \in C$. Then, as $\alpha \in [\mu, 1)$, $\mu \in [\max\{0, 1 - (\frac{\lambda q}{d_q})^{\frac{1}{q-1}}\}, 1)$, $T_\alpha$ is a nonexpansive mapping such that $F(T_\alpha) = F(T)$.*

*Proof* For any $x, y \in C$, we have, by Lemma 2.5 and (1.10),

$$\begin{aligned}
\|T_\alpha x - T_\alpha y\|^q &= \|\alpha x + (1-\alpha)Tx - \alpha y - (1-\alpha)Ty\|^q \\
&= \|x - y - (1-\alpha)[x - y - (Tx - Ty)]\|^q \\
&\leq \|x-y\|^q - q(1-\alpha)\langle (I-T)x - (I-T)y, j_q(x-y) \rangle \\
&\quad + d_q(1-\alpha)^q \|x - y - (Tx - Ty)\|^q \\
&\leq \|x-y\|^q - \lambda q(1-\alpha) \|x - y - (Tx - Ty)\|^q \\
&\quad + d_q(1-\alpha)^q \|x - y - (Tx - Ty)\|^q \\
&= \|x-y\|^q - (1-\alpha)[\lambda q - d_q(1-\alpha)^{q-1}] \|x - y - (Tx - Ty)\|^q \\
&\leq \|x-y\|^q,
\end{aligned}$$

which shows that $T_\alpha$ is a nonexpansive mapping.

It is clear that $x = T_\alpha x \Leftrightarrow x = Tx$. This proves our assertions. □

**Lemma 2.8** (Xu [18]) *Let $\{a_n\}$ be a sequence of nonnegative real numbers such that*

$$a_{n+1} \leq (1-\gamma_n)a_n + \delta_n, \quad n \geq 0,$$

*where $\{\gamma_n\}$ is a sequence in $(0,1)$ and $\{\delta_n\}$ is a sequence in $\mathbb{R}$ such that*

(i) $\lim_{n\to\infty} \gamma_n = 0$ and $\sum_{n=0}^\infty \gamma_n = \infty$;

(ii) $\limsup_{n\to\infty} \frac{\delta_n}{\gamma_n} \leq 0$ or $\sum_{n=1}^\infty |\delta_n| < \infty$.

*Then $\lim_{n\to\infty} a_n = 0$.*

**Lemma 2.9** *Let $E$ be a real $q$-uniformly smooth Banach space with constant $d_q$, $q > 1$, and let $C$ be a nonempty closed convex subset of $E$. Suppose that $T_i : C \to C$ are $k_i$-strict pseudocontractions for $k_i \in (0,1)$ ($i = 1, 2, \ldots, N$). Let $T_{\alpha_i} = \alpha_i I + (1-\alpha_i)T_i$, $k_i < \alpha_i < 1$ ($i = 1, 2, \ldots, N$). If $\bigcap_{i=1}^N F(T_i) \neq \emptyset$, then, as $\alpha_i \in [\mu, 1)$, $\mu \in [\max\{0, 1 - (\frac{\lambda q}{d_q})^{\frac{1}{q-1}}\}, 1)$, we have*

$$F(T_{\alpha_1} T_{\alpha_2} \cdots T_{\alpha_N}) = \bigcap_{i=1}^N F(T_{\alpha_i}).$$

*Proof* We prove it by induction. For $N = 2$, set $T_{\alpha_1} = \alpha_1 I + (1-\alpha_1)T_1$, $T_{\alpha_2} = \alpha_2 I + (1-\alpha_2)T_2$, $k_i < \alpha_i < 1$, $i = 1, 2$. Obviously,

$$F(T_{\alpha_1}) \cap F(T_{\alpha_2}) \subset F(T_{\alpha_1} T_{\alpha_2}).$$

Now we prove

$$F(T_{\alpha_1} T_{\alpha_2}) \subset F(T_{\alpha_1}) \cap F(T_{\alpha_2}).$$

For all $y \in F(T_{\alpha_1} T_{\alpha_2})$, $T_{\alpha_1} T_{\alpha_2} y = y$, if $T_{\alpha_2} y = y$, then $T_{\alpha_1} y = y$, the conclusion holds. In fact, we can claim that $T_{\alpha_2} y = y$. From Lemma 2.7, we know that $T_{\alpha_2}$ is nonexpansive and $F(T_{\alpha_1}) \cap F(T_{\alpha_2}) = F(T_1) \cap F(T_2) \neq \emptyset$.

Take $x \in F(T_{\alpha_1}) \cap F(T_{\alpha_2})$, then, by Lemma 2.5 and (1.10), we have

$$\begin{aligned}
\|x - y\|^q &= \|x - T_{\alpha_1} T_{\alpha_2} y\|^q \\
&= \|x - [\alpha_1(T_{\alpha_2} y) + (1 - \alpha_1) T_1 T_{\alpha_2} y]\|^q \\
&= \|x - T_{\alpha_2} y - (1 - \alpha_1)[x - T_{\alpha_2} y - (x - T_1 T_{\alpha_2} y)]\|^q \\
&\leq \|x - T_{\alpha_2} y\|^q - q(1 - \alpha_1)\langle x - T_{\alpha_2} y - (x - T_1 T_{\alpha_2} y), j_q(x - T_{\alpha_2} y)\rangle \\
&\quad + d_q(1 - \alpha_1)^q \|x - T_{\alpha_2} y - (x - T_1 T_{\alpha_2} y)\|^q \\
&= \|T_{\alpha_2} x - T_{\alpha_2} y\|^q - q(1 - \alpha_1)\langle (I - T_1)x - (I - T_1) T_{\alpha_2} y, j_q(x - T_{\alpha_2} y)\rangle \\
&\quad + d_q(1 - \alpha_1)^q \|x - T_{\alpha_2} y - (x - T_1 T_{\alpha_2} y)\|^q \\
&\leq \|T_{\alpha_2} x - T_{\alpha_2} y\|^q - \lambda q(1 - \alpha_1)\|x - T_{\alpha_2} y - (x - T_1 T_{\alpha_2} y)\|^q \\
&\quad + d_q(1 - \alpha_1)^q \|x - T_{\alpha_2} y - (x - T_1 T_{\alpha_2} y)\|^q \\
&\leq \|x - y\|^q - \lambda q(1 - \alpha_1)\|T_1 T_{\alpha_2} y - T_{\alpha_2} y\|^q \\
&\quad + d_q(1 - \alpha_1)^q \|T_1 T_{\alpha_2} y - T_{\alpha_2} y\|^q \\
&= \|x - y\|^q - (1 - \alpha_1)[\lambda q - d_q(1 - \alpha_1)^{q-1}]\|T_1 T_{\alpha_2} y - T_{\alpha_2} y\|^q.
\end{aligned}$$

So, we get

$$\|T_1 T_{\alpha_2} y - T_{\alpha_2} y\|^q \leq 0.$$

Namely $T_1 T_{\alpha_2} y = T_{\alpha_2} y$, that is,

$$T_{\alpha_2} y \in F(T_1) = F(T_{\alpha_1}), \qquad T_{\alpha_2} y = T_{\alpha_1} T_{\alpha_2} y = y.$$

Suppose that the conclusion holds for $N = k$, we prove that

$$F(T_{\alpha_1} T_{\alpha_2} \cdots T_{\alpha_{k+1}}) = \bigcap_{i=1}^{k+1} F(T_{\alpha_i}).$$

It suffices to verify

$$F(T_{\alpha_1} T_{\alpha_2} \cdots T_{\alpha_{k+1}}) \subset \bigcap_{i=1}^{k+1} F(T_{\alpha_i}).$$

For all $y \in F(T_{\alpha_1} T_{\alpha_2} \cdots T_{\alpha_{k+1}})$, $T_{\alpha_1} T_{\alpha_2} \cdots T_{\alpha_{k+1}} y = y$. Using Lemma 2.5 and (1.10) again, take $x \in \bigcap_{i=1}^{k+1} F(T_{\alpha_i})$, then

$$\begin{aligned}
\|x - y\|^q &= \|x - T_{\alpha_1} T_{\alpha_2} \cdots T_{\alpha_{k+1}} y\|^q \\
&= \|x - [\alpha_1(T_{\alpha_2} \cdots T_{\alpha_{k+1}} y) + (1 - \alpha_1) T_1 T_{\alpha_2} \cdots T_{\alpha_{k+1}} y]\|^q \\
&= \|x - T_{\alpha_2} \cdots T_{\alpha_{k+1}} y - (1 - \alpha_1)[x - T_{\alpha_2} \cdots T_{\alpha_{k+1}} y - (x - T_1 T_{\alpha_2} \cdots T_{\alpha_{k+1}} y)]\|^q \\
&\leq \|x - T_{\alpha_2} \cdots T_{\alpha_{k+1}} y\|^q - q(1 - \alpha_1)\langle x - T_{\alpha_2} \cdots T_{\alpha_{k+1}} y - (x - T_1 T_{\alpha_2} \cdots T_{\alpha_{k+1}} y),
\end{aligned}$$

$$j_q(x - T_{\alpha_2} \cdots T_{\alpha_{k+1}} y)\rangle$$
$$+ d_q(1-\alpha_1)^q \|x - T_{\alpha_2} \cdots T_{\alpha_{k+1}} y - (x - T_1 T_{\alpha_2} \cdots T_{\alpha_{k+1}} y)\|^q$$
$$= \|x - T_{\alpha_2} \cdots T_{\alpha_{k+1}} y\|^q - q(1-\alpha_1)\langle (I - T_1)x - (I - T_1)T_{\alpha_2} \cdots T_{\alpha_{k+1}} y,$$
$$j_q(x - T_{\alpha_2} \cdots T_{\alpha_{k+1}} y)\rangle$$
$$+ d_q(1-\alpha_1)^q \|x - T_{\alpha_2} \cdots T_{\alpha_{k+1}} y - (x - T_1 T_{\alpha_2} \cdots T_{\alpha_{k+1}} y)\|^q$$
$$\leq \|x - T_{\alpha_2} \cdots T_{\alpha_{k+1}} y\|^q$$
$$- \lambda q(1-\alpha_1) \|x - T_{\alpha_2} \cdots T_{\alpha_{k+1}} y - (x - T_1 T_{\alpha_2} \cdots T_{\alpha_{k+1}} y)\|^q$$
$$+ d_q(1-\alpha_1)^q \|x - T_{\alpha_2} \cdots T_{\alpha_{k+1}} y - (x - T_1 T_{\alpha_2} \cdots T_{\alpha_{k+1}} y)\|^q$$
$$\leq \|x - y\|^q - \lambda q(1-\alpha_1) \|T_1 T_{\alpha_2} \cdots T_{\alpha_{k+1}} y - T_{\alpha_2} \cdots T_{\alpha_{k+1}} y\|^q$$
$$+ d_q(1-\alpha_1)^q \|T_1 T_{\alpha_2} \cdots T_{\alpha_{k+1}} y - T_{\alpha_2} \cdots T_{\alpha_{k+1}} y\|^q$$
$$= \|x - y\|^q - (1-\alpha_1)\big[\lambda q - d_q(1-\alpha_1)^{q-1}\big]\|T_1 T_{\alpha_2} \cdots T_{\alpha_{k+1}} y - T_{\alpha_2} \cdots T_{\alpha_{k+1}} y\|^q.$$

So, we get

$$\|T_1 T_{\alpha_2} \cdots T_{\alpha_{k+1}} y - T_{\alpha_2} \cdots T_{\alpha_{k+1}} y\|^q \leq 0.$$

Thus, $T_1 T_{\alpha_2} \cdots T_{\alpha_{k+1}} y = T_{\alpha_2} \cdots T_{\alpha_{k+1}} y$, that is, $T_{\alpha_2} \cdots T_{\alpha_{k+1}} y \in F(T_1) = F(T_{\alpha_1})$. Namely,

$$T_{\alpha_2} \cdots T_{\alpha_{k+1}} y = T_{\alpha_1} T_{\alpha_2} \cdots T_{\alpha_{k+1}} y = y. \tag{2.1}$$

From (2.1) and inductive assumption, we get

$$y \in F(T_{\alpha_2} \cdots T_{\alpha_{k+1}}) = \bigcap_{i=2}^{k+1} F(T_{\alpha_i}),$$

that is,

$$T_{\alpha_i} y = y, \quad i = 2, \ldots, k+1.$$

Substituting it into (2.1), we obtain $T_{\alpha_1} T_{\alpha_i} y = y$, $i = 2, \ldots, k+1$, that is, $T_{\alpha_1} y = y$, $y \in F(T_{\alpha_1})$, and hence

$$y \in \bigcap_{i=1}^{k+1} F(T_{\alpha_i}). \qquad \square$$

**Lemma 2.10** (Ali *et al.* [9]) *Let E be a real q-uniformly smooth Banach space with constant $d_q$, $q > 1$. Let $f : E \to E$ be a contraction mapping with constant $\alpha \in (0,1)$. Let $T : E \to E$ be a nonexpansive mapping such that $F(T) \neq \emptyset$, and let $A : E \to E$ be an $\eta$-strongly accretive mapping which is also k-Lipschitzian. Let $\mu \in (0, \min\{1, (\frac{q\eta}{d_q k^q})^{\frac{1}{q-1}}\})$ and $\tau := \mu(\eta - \frac{\mu^{q-1} d_q k^q}{q})$. For each $t \in (0,1)$ and $\gamma \in (0, \frac{\tau}{\alpha})$, the path $\{x_t\}$ defined by*

$$x_t = t\gamma f(x_t) + (I - t\mu A)Tx_t$$

*converges strongly as* $t \to 0$ *to a fixed point* $x^*$ *of* $T$, *which solves the variational inequality*

$$\langle (\mu A - \gamma f)x^*, j_q(x^* - z) \rangle \leq 0, \quad \forall z \in F(T).$$

**Lemma 2.11** (Chang *et al.* [19]) *Let E be a real Banach space with a uniformly Gâteaux differentiable norm. Then the generalized duality mapping* $J_q : E \to 2^{E^*}$ *is single-valued and uniformly continuous on each bounded subset of E from the norm topology of E to the weak* topology of $E^*$.*

**Lemma 2.12** (Zhou *et al.* [20]) *Let* $\alpha$ *be a real number, and let a sequence* $\{a_n\} \in l^\infty$ *satisfy the condition* $\mu_n(a_n) \leq \alpha$ *for all Banach limit* $\mu$. *If* $\limsup_{n\to\infty}(a_{n+N} - a_n) \leq 0$, *then* $\limsup_{n\to\infty} a_n \leq \alpha$.

**Lemma 2.13** (Mitrinović [21]) *Suppose that* $q > 1$. *Then, for any arbitrary positive real numbers* $x, y$, *the following inequality holds*:

$$xy \leq \frac{1}{q}x^q + \left(\frac{q-1}{q}\right)x^{\frac{q}{q-1}}.$$

## 3 Synchronal algorithm

**Theorem 3.1** *Let E a real q-uniformly smooth Banach space, and let C be a nonempty closed convex subset of E. Let* $T_i : C \to C$ *be* $k_i$-*strict pseudocontractions for* $k_i \in (0, 1)$ ($i = 1, 2, \ldots, N$) *such that* $\bigcap_{i=1}^{N} F(T_i) \neq \emptyset$. *Let f be a contraction with coefficient* $\beta \in (0, 1)$, *and let* $\{\lambda_i\}_{i=1}^{N}$ *be a sequence of positive numbers such that* $\sum_{i=1}^{N} \lambda_i = 1$. *Let* $G : C \to C$ *be an* $\eta$-*strongly accretive and L-Lipschitzian operator with* $L > 0$, $\eta > 0$. *Assume that* $0 < \mu < (q\eta/d_q L^q)^{1/q-1}$, $0 < \gamma < \mu(\eta - d_q \mu^{q-1} L^q/q)/\beta = \tau/\beta$. *Let* $\{\alpha_n\}$ *and* $\{\beta_n\}$ *be sequences in* $(0, 1)$ *satisfying the following conditions*:

(K1) $\lim_{n\to\infty} \alpha_n = 0$, $\sum_{n=0}^{\infty} \alpha_n = \infty$;

(K2) $\sum_{n=1}^{\infty} |\alpha_{n+1} - \alpha_n| < \infty$, $\sum_{n=1}^{\infty} |\beta_{n+1} - \beta_n| < \infty$;

(K3) $0 < k \leq \beta_n < a < 1$, *where* $k = \min\{k_i : 1 \leq i \leq N\}$;

(K4) $\alpha_n, \beta_n \in [\mu, 1)$, *where* $\mu \in [\max\{0, 1 - (\frac{\lambda q}{d_q})^{\frac{1}{q-1}}\}, 1)$.

*Let* $\{x_n\}$ *be a sequence defined by algorithm* (1.17), *then* $\{x_n\}$ *converges strongly to a common fixed point of* $\{T_i\}_{i=1}^{N}$, *which solves the variational inequality* (1.16).

*Proof* Let $T := \sum_{i=1}^{N} \lambda_i T_i$, then by Lemma 2.1 we conclude that $T$ is a $k$-strict pseudocontraction and $F(T) = \bigcap_{i=1}^{N} F(T_i)$. We can rewrite algorithm (1.17) as follows:

$$\begin{cases} x_0 = x \in C \quad \text{chosen arbitrarily}, \\ T^{\beta_n} = \beta_n I + (1 - \beta_n)T, \\ x_{n+1} = \alpha_n \gamma f(x_n) + (I - \alpha_n \mu G)T^{\beta_n} x_n, \quad n \geq 0. \end{cases}$$

Furthermore, by Lemma 2.7 we have that $T^{\beta_n}$ is a nonexpansive mapping and $F(T^{\beta_n}) = F(T)$. From condition (K1) we may assume, without loss of generality, that $\alpha_n \in (0, \min\{1, \frac{1}{\tau}\})$. Let $p \in \bigcap_{i=1}^{N} F(T_i)$, then the sequence $\{x_n\}$ satisfies

$$\|x_n - p\| \leq \max\left\{\|x_0 - p\|, \frac{\|\gamma f(p) - \mu G p\|}{\tau - \gamma \beta}\right\}, \quad \forall n \geq 0.$$

We prove this by mathematical induction as follows.

Obviously, it is true for $n = 0$. Assume that it is true for $n = k$ for some $k \in \mathbb{N}$.

From (1.17) and Lemma 2.6, we have

$$\begin{aligned}
\|x_{k+1} - p\| &= \left\|\alpha_k \gamma f(x_k) + (I - \alpha_k \mu G) T^{\beta_k} x_k - p\right\| \\
&= \left\|\alpha_k [\gamma f(x_k) - \mu G p] + (I - \alpha_k \mu G) T^{\beta_k} x_k - (I - \alpha_k \mu G) p\right\| \\
&\leq (1 - \alpha_k \tau)\|x_k - p\| + \alpha_k \left\|\gamma [f(x_k) - f(p)] + \gamma f(p) - \mu G p\right\| \\
&\leq (1 - \alpha_k \tau)\|x_k - p\| + \alpha_k \gamma \beta \|x_k - p\| + \alpha_k \|\gamma f(p) - \mu G p\| \\
&= [1 - \alpha_k (\tau - \gamma \beta)]\|x_k - p\| + \alpha_k (\tau - \gamma \beta) \frac{\|\gamma f(p) - \mu G p\|}{\tau - \gamma \beta} \\
&\leq \max \left\{ \|x_k - p\|, \frac{\|\gamma f(p) - \mu G p\|}{\tau - \gamma \beta} \right\}.
\end{aligned}$$

Hence the proof. Thus, the sequence $\{x_n\}$ is bounded and so are $\{Tx_n\}$, $\{GT^{\beta_n} x_n\}$ and $\{f(x_n)\}$.

Observe that

$$\begin{aligned}
x_{n+2} - x_{n+1} &= [\alpha_{n+1} \gamma f(x_{n+1}) + (I - \alpha_{n+1} \mu G) T^{\beta_{n+1}} x_{n+1}] \\
&\quad - [\alpha_n \gamma f(x_n) + (I - \alpha_n \mu G) T^{\beta_n} x_n] \\
&= [\alpha_{n+1} \gamma f(x_{n+1}) - \alpha_{n+1} \gamma f(x_n)] + [\alpha_{n+1} \gamma f(x_n) - \alpha_n \gamma f(x_n)] \\
&\quad + [(I - \alpha_{n+1} \mu G) T^{\beta_{n+1}} x_{n+1} - (I - \alpha_{n+1} \mu G) T^{\beta_n} x_n] \\
&\quad + [\alpha_n \mu G T^{\beta_n} x_n - \alpha_{n+1} \mu G T^{\beta_n} x_n] \\
&= \alpha_{n+1} \gamma [f(x_{n+1}) - f(x_n)] + [(I - \alpha_{n+1} \mu G) T^{\beta_{n+1}} x_{n+1} \\
&\quad - (I - \alpha_{n+1} \mu G) T^{\beta_n} x_n] + (\alpha_{n+1} - \alpha_n) \gamma f(x_n) \\
&\quad + (\alpha_n - \alpha_{n+1}) \mu G T^{\beta_n} x_n,
\end{aligned}$$

so that

$$\begin{aligned}
\|x_{n+2} - x_{n+1}\| &\leq \alpha_{n+1} \gamma \beta \|x_{n+1} - x_n\| + (1 - \alpha_{n+1} \tau) \left\|T^{\beta_{n+1}} x_{n+1} - T^{\beta_n} x_n\right\| \\
&\quad + |\alpha_{n+1} - \alpha_n|\left(\gamma \|f(x_n)\| + \mu \|GT^{\beta_n} x_n\|\right) \\
&\leq \alpha_{n+1} \gamma \beta \|x_{n+1} - x_n\| + (1 - \alpha_{n+1} \tau) \left\|T^{\beta_{n+1}} x_{n+1} - T^{\beta_n} x_n\right\| \\
&\quad + |\alpha_{n+1} - \alpha_n| M_1,
\end{aligned} \qquad (3.1)$$

where $M_1$ is an appropriate constant such that $M_1 \geq \sup_{n \geq 1}\{\gamma \|f(x_n)\| + \mu \|GT^{\beta_n} x_n\|\}$.

On the other hand, we note that

$$\begin{aligned}
\left\|T^{\beta_{n+1}} x_{n+1} - T^{\beta_n} x_n\right\| &\leq \left\|T^{\beta_{n+1}} x_{n+1} - T^{\beta_{n+1}} x_n\right\| + \left\|T^{\beta_{n+1}} x_n - T^{\beta_n} x_n\right\| \\
&\leq \|x_{n+1} - x_n\| + \left\|[\beta_{n+1} x_n + (1 - \beta_{n+1}) T x_n] \right. \\
&\quad \left. - [\beta_n x_n + (1 - \beta_n) T x_n]\right\|
\end{aligned}$$

$$= \|x_n - x_n\| + \|\beta_{n+1}(x_n - Tx_n) - \beta_n(x_n - Tx_n)\|$$
$$\leq \|x_{n+1} - x_n\| + |\beta_{n+1} - \beta_n|M_2, \qquad (3.2)$$

where $M_2$ is an appropriate constant such that $M_2 \geq \sup_{n \geq 1}\{\|x_n - Tx_n\|\}$.

Now, substituting (3.2) into (3.1) yields

$$\|x_{n+2} - x_{n+1}\| \leq \alpha_{n+1}\gamma\beta\|x_{n+1} - x_n\| + (1 - \alpha_{n+1}\tau)\|x_{n+1} - x_n\| + |\alpha_{n+1} - \alpha_n|M_1$$
$$+ |\beta_{n+1} - \beta_n|M_2$$
$$\leq [1 - \alpha_{n+1}(\tau - \gamma\beta)]\|x_{n+1} - x_n\| + (|\alpha_{n+1} - \alpha_n| + |\beta_{n+1} - \beta_n|)M_3,$$

where $M_3$ is an appropriate constant such that $M_3 \geq \max\{M_1, M_2\}$.

By Lemma 2.8 and conditions (K1), (K2), we have

$$\|x_{n+1} - x_n\| \to 0 \quad \text{as } n \to \infty. \qquad (3.3)$$

From (1.17) and condition (K1), we have

$$\|x_{n+1} - T^{\beta_n}x_n\| = \|\alpha_n\gamma f(x_n) + (I - \alpha_n\mu G)T^{\beta_n}x_n - T^{\beta_n}x_n\|$$
$$\leq \alpha_n\|\gamma f(x_n) + \mu G T^{\beta_n}x_n\| \to 0 \quad \text{as } n \to \infty. \qquad (3.4)$$

On the other hand,

$$\|x_{n+1} - T^{\beta_n}x_n\| = \|x_{n+1} - [\beta_n x_n + (1 - \beta_n)Tx_n]\|$$
$$= \|(x_{n+1} - x_n) + (1 - \beta_n)(x_n - Tx_n)\|$$
$$\geq (1 - \beta_n)\|x_n - Tx_n\| - \|x_{n+1} - x_n\|,$$

which implies, by condition (K3), that

$$\|x_n - Tx_n\| \leq \frac{1}{1 - \beta_n}(\|x_{n+1} - x_n\| + \|x_{n+1} - T^{\beta_n}x_n\|)$$
$$\leq \frac{1}{1 - a}(\|x_{n+1} - x_n\| + \|x_{n+1} - T^{\beta_n}x_n\|).$$

Hence, from (3.3) and (3.4), we have

$$\|x_n - Tx_n\| \to 0 \quad \text{as } n \to \infty. \qquad (3.5)$$

From the boundedness of $\{x_n\}$, without loss of generality, we may assume that $x_n \rightharpoonup p$. Hence, by Lemma 2.2 and (3.5), we obtain $Tp = p$. So, we have

$$\omega_\omega(x_n) \subset F(T). \qquad (3.6)$$

We now prove that $\limsup_{n\to\infty} \langle(\gamma f - \mu G)x^*, j_q(x_{n+1} - x^*)\rangle \leq 0$, where $x^*$ is obtained in Lemma 2.10. Put $a_n := \langle(\gamma f - \mu G)x^*, j_q(x_n - x^*)\rangle$, $n \geq 1$.

Define a map $\phi : E \to \mathbb{R}$ by

$$\phi(x) = \mu_n \|x_n - x\|^q, \quad \forall x \in E.$$

Then $\phi$ is continuous, convex, and $\phi(x) \to \infty$ as $\|x\| \to \infty$. Since $E$ is reflexive, there exists $y^* \in C$ such that $\phi(y_*) = \min_{z \in C} \phi(z)$. Hence the set

$$K^* := \left\{ y^* \in C : \phi(y^*) = \min_{z \in C} \phi(z) \right\} \neq \emptyset.$$

Therefore, applying Lemma 2.4, we have $K^* \cap F(T^{\beta_n}) \neq \emptyset$. Without loss of generality, assume $x^* = y^* \in K^* \cap F(T^{\beta_n})$. Let $t \in (0,1)$. Then it follows that $\phi(x^*) \leq \phi(x^* + t(\gamma f - \mu G)x^*)$, and using Lemma 2.3, we obtain that

$$\|x_n - x^* - t(\gamma f - \mu G)x^*\|^q \leq \|x_n - x^*\|^q - qt\langle (\gamma f - \mu G)x^*, j_q(x_n - x^* - t(\gamma f - \mu G)x^*) \rangle.$$

This implies that

$$\mu_n \langle (\gamma f - \mu G)x^*, j_q(x_n - x^* - t(\gamma f - \mu G)x^*) \rangle \leq 0.$$

By Lemma 2.11, $j_q$ is norm-to-weak* uniformly continuous on a bounded subset of $E$, so we obtain, as $t \to 0$, that

$$\langle (\gamma f - \mu G)x^*, j_q(x_n - x^*) \rangle - \langle (\gamma f - \mu G)x^*, j_q(x_n - x^* - t(\gamma f - \mu G)x^*) \rangle \to 0.$$

Hence, for all $\epsilon > 0$, there exists $\delta_\epsilon > 0$ such that $\forall t \in (0, \delta_\epsilon)$ and for all $n \geq 1$,

$$\langle (\gamma f - \mu G)x^*, j_q(x_n - x^*) \rangle - \langle (\gamma f - \mu G)x^*, j_q(x_n - x^* - t(\gamma f - \mu G)x^*) \rangle < \epsilon.$$

Consequently,

$$\mu_n \langle (\gamma f - \mu G)x^*, j_q(x_n - x^*) \rangle \leq \mu_n \langle (\gamma f - \mu G)x^*, j_q(x_n - x^* - t(\gamma f - \mu G)x^*) \rangle + \epsilon \leq \epsilon.$$

Since $\epsilon$ is arbitrary, we have

$$\mu_n \langle (\gamma f - \mu G)x^*, j_q(x_n - x^*) \rangle \leq 0.$$

Thus, $\mu_n(a_n) \leq 0$ for any Banach limit $\mu$.

Furthermore, by (3.3), $\|x_{n+1} - x_n\| \to 0$ as $n \to \infty$. We therefore conclude that

$$\limsup_{n \to \infty}(a_{n+1} - a_n) = \limsup_{n \to \infty}\left(\langle (\gamma f - \mu G)x^*, j_q(x_{n+1} - x^*) \rangle - \langle (\gamma f - \mu G)x^*, j_q(x_n - x^*) \rangle\right)$$
$$= \limsup_{n \to \infty}\left(\langle (\gamma f - \mu G)x^*, j_q(x_{n+1} - x^*) - j_q(x_n - x^*) \rangle\right) = 0.$$

Hence, by Lemma 2.12 we obtain $\limsup_{n \to \infty} a_n \leq 0$, that is,

$$\limsup_{n \to \infty}\langle (\gamma f - \mu G)x^*, j_q(x_n - x^*) \rangle \leq 0. \tag{3.7}$$

From (1.17), Lemmas 2.3, 2.6 and 2.13, we have

$$\begin{aligned}
\|x_{n+1} - x^*\|^q &= \langle x_{n+1} - x^*, j_q(x_{n+1} - x^*) \rangle \\
&= \langle \alpha_n[\gamma f(x_n) - \mu G x^*] + (I - \alpha_n \mu G)(T^{\beta_n} x_n - x^*), j_q(x_{n+1} - x^*) \rangle \\
&= \alpha_n \langle \gamma f(x_n) - \mu G x^*, j_q(x_{n+1} - x^*) \rangle \\
&\quad + \langle (I - \alpha_n \mu G)(T^{\beta_n} x_n - x^*), j_q(x_{n+1} - x^*) \rangle \\
&= \alpha_n \langle \gamma f(x_n) - \gamma f(x^*), j_q(x_{n+1} - x^*) \rangle \\
&\quad + \alpha_n \langle \gamma f(x^*) - \mu G x^*, j_q(x_{n+1} - x^*) \rangle \\
&\quad + \langle (I - \alpha_n \mu G)(T^{\beta_n} x_n - x^*), j_q(x_{n+1} - x^*) \rangle \\
&\leq \alpha_n \gamma \|f(x_n) - f(x^*)\| \|x_{n+1} - x^*\|^{q-1} + \alpha_n \langle (\gamma f - \mu G) x^*, j_q(x_{n+1} - x^*) \rangle \\
&\quad + \|(I - \alpha_n \mu G)T^{\beta_n} x_n - (I - \alpha_n \mu G)x^*\| \|x_{n+1} - x^*\|^{q-1} \\
&\leq (1 - \alpha_n \tau)\|x_n - x^*\| \|x_{n+1} - x^*\|^{q-1} + \alpha_n \gamma \beta \|x_n - x^*\| \|x_{n+1} - x^*\|^{q-1} \\
&\quad + \alpha_n \langle (\gamma f - \mu G)x^*, j_q(x_{n+1} - x^*) \rangle \\
&= [1 - \alpha_n(\tau - \gamma\beta)]\|x_n - x^*\| \|x_{n+1} - x^*\|^{q-1} \\
&\quad + \alpha_n \langle (\gamma f - \mu G)x^*, j_q(x_{n+1} - x^*) \rangle \\
&\leq [1 - \alpha_n(\tau - \gamma\beta)]\left[\frac{1}{q}\|x_n - x^*\|^q + \left(\frac{q-1}{q}\right)\|x_{n+1} - x^*\|^q\right] \\
&\quad + \alpha_n \langle (\gamma f - \mu G)x^*, j_q(x_{n+1} - x^*) \rangle.
\end{aligned}$$

This implies that

$$\begin{aligned}
\|x_{n+1} - x^*\|^q &\leq \frac{1 - \alpha_n(\tau - \gamma\beta)}{1 + \alpha_n(q-1)(\tau - \gamma\beta)} \|x_n - x^*\|^q \\
&\quad + \frac{q\alpha_n}{1 + \alpha_n(q-1)(\tau - \gamma\beta)} \langle (\gamma f - \mu G)x^*, j_q(x_{n+1} - x^*) \rangle \\
&\leq [1 - \alpha_n(\tau - \gamma\beta)]\|x_n - x^*\|^q \\
&\quad + \frac{q\alpha_n}{1 + \alpha_n(q-1)(\tau - \gamma\beta)} \langle (\gamma f - \mu G)x^*, j_q(x_{n+1} - x^*) \rangle \\
&\leq (1 - \gamma_n)\|x_n - x^*\|^q + \delta_n,
\end{aligned}$$

where $\gamma_n := \alpha_n(\tau - \gamma\beta)$ and $\delta_n := \frac{q\alpha_n}{1+\alpha_n(q-1)(\tau-\gamma\beta)} \langle (\gamma f - \mu G)x^*, j_q(x_{n+1} - x^*) \rangle$. From (K1), $\lim_{n \to \infty} \gamma_n = 0$, $\sum_{n=0}^{\infty} \gamma_n = \infty$. $\frac{\delta_n}{\gamma_n} = \frac{q}{[1+\alpha_n(q-1)(\tau-\gamma\beta)](\tau-\gamma\beta)} \langle (\gamma f - \mu G)x^*, j_q(x_{n+1} - x^*) \rangle$. So, $\limsup_{n \to \infty} \frac{\delta_n}{\gamma_n} \leq 0$. Hence, by Lemma 2.8, we have that $x_n \to x^*$ as $n \to \infty$. This completes the proof. □

## 4 Cyclic algorithm

**Theorem 4.1** *Let $E$ be a real $q$-uniformly smooth Banach space, and let $C$ be a nonempty closed convex subset of $E$. Let $T_i : C \to C$ be $k_i$-strict pseudocontractions for $k_i \in (0,1)$ $(i = 1, 2, \ldots, N)$ such that $\bigcap_{i=1}^{N} F(T_i) \neq \emptyset$, let $f$ be a contraction with coefficient $\beta \in (0,1)$. Let $G : C \to C$ be an $\eta$-strongly accretive and $L$-Lipschitzian operator with $L > 0$, $\eta > 0$.*

*Assume that $0 < \mu < (q\eta/d_q L^q)^{1/q-1}$, $0 < \gamma < \mu(\eta - d_q \mu^{q-1} L^q/q)/\beta = \tau/\beta$. Let $\{\alpha_n\}$ and $\{\beta_n\}$ be sequences in $(0,1)$ satisfying the following conditions*:

(K1′) $\lim_{n\to\infty} \alpha_n = 0$, $\sum_{n=0}^{\infty} \alpha_n = \infty$;
(K2′) $\sum_{n=1}^{\infty} |\alpha_{n+1} - \alpha_n| < \infty$ or $\lim_{n\to\infty} \frac{\alpha_n}{\alpha_{n+N}} = 1$;
(K3′) $\beta_{[n]} \in [k, 1)$, $\forall n \geq 0$, where $k = \min\{k_i : 1 \leq i \leq N\}$;
(K4′) $\alpha_n, \beta_n \in [\mu, 1)$, where $\mu \in [\max\{0, 1 - (\frac{\lambda q}{d_q})^{\frac{1}{q-1}}\}, 1)$.

*Let $\{x_n\}$ be a sequence defined by algorithm* (1.18), *then $\{x_n\}$ converges strongly to a common fixed point of $\{T_i\}_{i=1}^{N}$, which solves the variational inequality* (1.16).

*Proof* From condition (K1′) we may assume, without loss of generality, that $\alpha_n \in (0, \min\{1, \frac{1}{\tau}\})$. Let $p \in \bigcap_{i=1}^{N} F(T_i)$, then the sequence $\{x_n\}$ satisfies

$$\|x_n - p\| \leq \max\left\{\|x_0 - p\|, \frac{\|\gamma f(p) - \mu G p\|}{\tau - \gamma \beta}\right\}, \quad \forall n \geq 0.$$

We prove this by mathematical induction as follows.

Obviously, it is true for $n = 0$. Assume it is true for $n = k$ for some $k \in \mathbb{N}$.
From (1.18) and Lemma 2.6, we have

$$\begin{aligned}
\|x_{k+1} - p\| &= \|\alpha_k \gamma f(x_k) + (I - \alpha_k \mu G) A_{[k+1]} x_k - p\| \\
&= \|\alpha_k [\gamma f(x_k) - \mu G p] + (I - \alpha_k \mu G) A_{[k+1]} x_k - (I - \alpha_k \mu G) p\| \\
&\leq (1 - \alpha_k \tau)\|x_k - p\| + \alpha_k \|\gamma [f(x_k) - f(p)] + \gamma f(p) - \mu G p\| \\
&\leq (1 - \alpha_k \tau)\|x_k - p\| + \alpha_k \gamma \beta \|x_k - p\| + \alpha_k \|\gamma f(p) - \mu G p\| \\
&= [1 - \alpha_k (\tau - \gamma \beta)]\|x_k - p\| + \alpha_k (\tau - \gamma \beta) \frac{\|\gamma f(p) - \mu G p\|}{\tau - \gamma \beta} \\
&\leq \max\left\{\|x_k - p\|, \frac{\|\gamma f(p) - \mu G p\|}{\tau - \gamma \beta}\right\}.
\end{aligned}$$

Hence the proof. Thus, the sequence $\{x_n\}$ is bounded and so are $\{T_{[n]} x_n\}$, $\{G A_{[n]} x_n\}$, $\{f(x_n)\}$, and $\{A_{[n]} x_n\}$.

From (1.18) and Lemma 2.6, we have

$$\begin{aligned}
\|x_{n+N+1} - x_{n+1}\| &= \|[\alpha_{n+N} \gamma f(x_{n+N}) + (I - \alpha_{n+N} \mu G) A_{[n+1]} x_{n+N}] \\
&\quad - [\alpha_n \gamma f(x_n) - (I - \alpha_n \mu G) A_{[n+1]} x_n]\| \\
&= \|\alpha_{n+N} \gamma f(x_{n+N}) - \alpha_{n+N} \gamma f(x_n) + \alpha_{n+N} \gamma f(x_n) \\
&\quad - \alpha_n \gamma f(x_n) + (I - \alpha_{n+N} \mu G) A_{[n+1]} x_{n+N} \\
&\quad - (I - \alpha_{n+N} \mu G) A_{[n+1]} x_n + (I - \alpha_{n+N} \mu G) A_{[n+1]} x_n \\
&\quad - (I - \alpha_n \mu G) A_{[n+1]} x_n\| \\
&= \|\alpha_{n+N} \gamma [f(x_{n+N}) - f(x_n)] + (\alpha_{n+N} - \alpha_n) \gamma f(x_n) \\
&\quad + (I - \alpha_{n+N} \mu G) A_{[n+1]} x_{n+N} - (I - \alpha_{n+N} \mu G) A_{[n+1]} x_n \\
&\quad + (\alpha_n - \alpha_{n+N}) \mu G A_{[n+1]} x_n\|
\end{aligned}$$

$$\leq \alpha_{n+N}\gamma\beta\|x_{n+N} - x_n\| + |\alpha_{n+N} - \alpha_n|\gamma\|f(x_n)\|$$

$$+ (1 - \alpha_{n+N}\tau)\|x_{n+N} - x_n\| + |\alpha_{n+N} - \alpha_n|\mu\|GA_{[n+1]}x_n\|$$

$$\leq [1 - \alpha_{n+N}(\tau - \gamma\beta)]\|x_{n+N} - x_n\| + |\alpha_{n+N} - \alpha_n|M_4,$$

where $M_4$ is an appropriate constant such that $M_4 \geq \sup_{n\geq 1}\{\mu\|GA_{[n+1]}x_n\| + \gamma\|f(x_n)\|\}$. By conditions (K1′), (K2′) and Lemma 2.8, we have

$$\|x_{n+N} - x_n\| \to 0 \quad \text{as } n \to \infty. \tag{4.1}$$

From (1.18) and condition (K1′), we have

$$\|x_{n+1} - A_{[n+1]}x_n\| = \|\alpha_n\gamma f(x_n) + (I - \alpha_n\mu G)A_{[n+1]}x_n - A_{[n+1]}x_n\|$$

$$= \alpha_n\|\gamma f(x_n) - \mu GA_{[n+1]}x_n\| \to 0 \quad \text{as } n \to \infty. \tag{4.2}$$

Recursively,

$$\|x_{n+N} - A_{[n+N]}x_{n+N-1}\| \to 0 \quad \text{as } n \to \infty,$$

$$\|x_{n+N-1} - A_{[n+N-1]}x_{n+N-2}\| \to 0 \quad \text{as } n \to \infty.$$

By condition (K3′) and Lemma 2.7, we know that $A_{[n+N]}$ is nonexpansive, so we get

$$\|A_{[n+N]}x_{n+N-1} - A_{[n+N]}A_{[n+N-1]}x_{n+N-2}\| \to 0 \quad \text{as } n \to \infty.$$

Proceeding accordingly, we have

$$\|A_{[n+N]}A_{[n+N-1]}x_{n+N-2} - A_{[n+N]}A_{[n+N-1]}A_{[n+N-2]}x_{n+N-3}\| \to 0 \quad \text{as } n \to \infty,$$

$$\vdots$$

$$\|A_{[n+N]}\cdots A_{[n+2]}x_{n+1} - A_{[n+N]}\cdots A_{[n+1]}x_n\| \to 0 \quad \text{as } n \to \infty.$$

Note that

$$\|x_{n+N} - A_{[n+N]}\cdots A_{[n+1]}x_n\| \leq \|x_{n+N} - A_{[n+N]}x_{n+N-1}\|$$

$$+ \|A_{[n+N]}x_{n+N-1} - A_{[n+N]}A_{[n+N-1]}x_{n+N-2}\|$$

$$+ \cdots$$

$$+ \|A_{[n+N]}\cdots A_{[n+2]}x_{n+1} - A_{[n+N]}\cdots A_{[n+1]}x_n\|.$$

From the above inequality, we obtain

$$\|x_{n+N} - A_{[n+N]}\cdots A_{[n+1]}x_n\| \to 0 \quad \text{as } n \to \infty.$$

Since

$$\|x_n - A_{[n+N]}\cdots A_{[n+1]}x_n\| \leq \|x_n - x_{n+N}\| + \|x_{n+N} - A_{[n+N]}\cdots A_{[n+1]}x_n\| \to 0 \quad \text{as } n \to \infty,$$

we conclude that

$$\|x_n - A_{[n+N]} \cdots A_{[n+1]} x_n\| \to 0 \quad \text{as } n \to \infty. \tag{4.3}$$

Take a subsequence $\{x_{n_j}\} \subset \{x_n\}$, by (4.3) we get

$$\|x_{n_j} - A_{[n_j+N]} \cdots A_{[n_j+1]} x_{n_j}\| \to 0 \quad \text{as } j \to \infty.$$

Notice that for each $n_j$, $A_{[n_j+N]} A_{[n_j+N-1]} \cdots A_{[n_j+1]}$ is some permutation of the mappings $A_1 A_2 \cdots A_N$. Since $A_1, A_2, \ldots, A_N$ are finite, all the finite permutations are $N!$, there must be some permutation appearing infinitely many times. Without loss of generality, suppose this permutation is $A_1 A_2 \cdots A_N$, we can take a subsequence $\{x_{n_{j_k}}\} \subset \{x_{n_j}\}$ such that $x_{n_{j_k}} \rightharpoonup q$ $(k \to \infty)$ and

$$\|x_{n_{j_k}} - A_1 A_2 \cdots A_N x_{n_{j_k}}\| \to 0 \quad \text{as } k \to \infty.$$

By Lemma 2.7, we conclude that $A_1, A_2, \ldots, A_N$ are all nonexpansive. It is clear that $A_{[n_j+N]} A_{[n_j+N-1]} \cdots A_{[n_j+1]}$ is nonexpansive, so is $A_1 A_2 \cdots A_N$. By Lemma 2.2, we have $A_1 A_2 \cdots A_N q = q$. From Lemmas 2.7 and 2.9, we obtain

$$q \in F(A_1 A_2 \cdots A_N) = \bigcap_{i=1}^{N} F(A_i) = \bigcap_{i=1}^{N} F(T_i),$$

that is,

$$\omega_\omega(x_n) \subset \bigcap_{i=1}^{N} F(T_i). \tag{4.4}$$

We now prove that $\limsup_{n \to \infty} \langle (\gamma f - \mu G) x^*, j_q(x_{n+1} - x^*) \rangle \leq 0$, where $x^*$ is obtained in Lemma 2.10. Put $a_n := \langle (\gamma f - \mu G) x^*, j_q(x_n - x^*) \rangle$, $n \geq 1$.

Define a map $\phi : E \to \mathbb{R}$ by

$$\phi(x) = \mu_n \|x_n - x\|^q, \quad \forall x \in E.$$

Then $\phi$ is continuous, convex, and $\phi(x) \to \infty$ as $\|x\| \to \infty$. Since $E$ is reflexive, there exists $y^* \in C$ such that $\phi(y^*) = \min_{z \in C} \phi(z)$. Hence the set

$$K^* := \left\{ y^* \in C : \phi(y^*) = \min_{z \in C} \phi(z) \right\} \neq \emptyset.$$

Therefore, applying Lemma 2.4, we have $K^* \cap F(T^{\beta_n}) \neq \emptyset$. Without loss of generality, assume $x^* = y^* \in K^* \cap F(T^{\beta_n})$. Let $t \in (0,1)$. Then it follows that $\phi(x^*) \leq \phi(x^* + t(\gamma f - \mu G)x^*)$, and using Lemma 2.3, we obtain that

$$\|x_n - x^* - t(\gamma f - \mu G)x^*\|^q \leq \|x_n - x^*\|^q - qt\langle (\gamma f - \mu G)x^*, j_q(x_n - x^* - t(\gamma f - \mu G)x^*) \rangle.$$

This implies that

$$\mu_n \langle (\gamma f - \mu G)x^*, j_q(x_n - x^* - t(\gamma f - \mu G)x^*) \rangle \leq 0.$$

By Lemma 2.11, $j_q$ is norm-to-weak* uniformly continuous on a bounded subset of $E$, so we obtain, as $t \to 0$, that

$$\langle (\gamma f - \mu G)x^*, j_q(x_n - x^*) \rangle - \langle (\gamma f - \mu G)x^*, j_q(x_n - x^* - t(\gamma f - \mu G)x^*) \rangle \to 0.$$

Hence, for all $\epsilon > 0$, there exists $\delta_\epsilon > 0$ such that $\forall t \in (0, \delta_\epsilon)$, and for all $n \geq 1$,

$$\langle (\gamma f - \mu G)x^*, j_q(x_n - x^*) \rangle - \langle (\gamma f - \mu G)x^*, j_q(x_n - x^* - t(\gamma f - \mu G)x^*) \rangle < \epsilon.$$

Consequently,

$$\mu_n \langle (\gamma f - \mu G)x^*, j_q(x_n - x^*) \rangle \leq \mu_n \langle (\gamma f - \mu G)x^*, j_q(x_n - x^* - t(\gamma f - \mu G)x^*) \rangle + \epsilon \leq \epsilon.$$

Since $\epsilon$ is arbitrary, we have

$$\mu_n \langle (\gamma f - \mu G)x^*, j_q(x_n - x^*) \rangle \leq 0.$$

Thus, $\mu_n(a_n) \leq 0$ for any Banach limit $\mu$.

Furthermore, by (4.1) $\|x_{n+N} - x_n\| \to 0$ as $n \to \infty$, we therefore conclude that

$$\limsup_{n \to \infty}(a_{n+N} - a_n) = \limsup_{n \to \infty} \left( \langle (\gamma f - \mu G)x^*, j_q(x_{n+N} - x^*) \rangle - \langle (\gamma f - \mu G)x^*, j_q(x_n - x^*) \rangle \right)$$

$$= \limsup_{n \to \infty} \left( \langle (\gamma f - \mu G)x^*, j_q(x_{n+N} - x^*) - j_q(x_n - x^*) \rangle \right) = 0.$$

Hence, by Lemma 2.12 we obtain $\limsup_{n \to \infty} a_n \leq 0$, that is,

$$\limsup_{n \to \infty} \langle (\gamma f - \mu G)x^*, j_q(x_n - x^*) \rangle \leq 0. \tag{4.5}$$

From (1.18), Lemmas 2.3, 2.6 and 2.13, we have

$$\|x_{n+1} - x^*\|^q = \langle x_{n+1} - x^*, j_q(x_{n+1} - x^*) \rangle$$

$$= \langle \alpha_n [\gamma f(x_n) - \mu G x^*] + (I - \alpha_n \mu G)(A_{[n+1]} x_n - x^*), j_q(x_{n+1} - x^*) \rangle$$

$$= \alpha_n \langle \gamma f(x_n) - \mu G x^*, j_q(x_{n+1} - x^*) \rangle$$

$$\quad + \langle (I - \alpha_n \mu G)(A_{[n+1]} x_n - x^*), j_q(x_{n+1} - x^*) \rangle$$

$$= \alpha_n \langle \gamma f(x_n) - \gamma f(x^*), j_q(x_{n+1} - x^*) \rangle + \alpha_n \langle \gamma f(x^*) - \mu G x^*, j_q(x_{n+1} - x^*) \rangle$$

$$\quad + \langle (I - \alpha_n \mu G)(A_{[n+1]} x_n - x^*), j_q(x_{n+1} - x^*) \rangle$$

$$\leq \alpha_n \gamma \|f(x_n) - f(x^*)\| \|x_{n+1} - x^*\|^{q-1} + \alpha_n \langle \gamma f - \mu G)x^*, j_q(x_{n+1} - x^*) \rangle$$

$$\quad + \|(I - \alpha_n \mu G)A_{[n+1]} x_n - (I - \alpha_n \mu G)x^*\| \|x_{n+1} - x^*\|^{q-1}$$

$$\leq (1 - \alpha_n \tau) \|x_n - x^*\| \|x_{n+1} - x^*\|^{q-1} + \alpha_n \gamma \beta \|x_n - x^*\| \|x_{n+1} - x^*\|^{q-1}$$

$$\quad + \alpha_n \langle (\gamma f - \mu G)x^*, j_q(x_{n+1} - x^*) \rangle$$

$$= \left[ 1 - \alpha_n(\tau - \gamma \beta) \right] \|x_n - x^*\| \|x_{n+1} - x^*\|^{q-1}$$

$$+ \alpha_n \langle (\gamma f - \mu G) x^*, j_q(x_{n+1} - x^*) \rangle$$
$$\leq [1 - \alpha_n(\tau - \gamma\beta)] \left[ \frac{1}{q} \|x_n - x^*\|^q + \left(\frac{q-1}{q}\right) \|x_{n+1} - x^*\|^q \right]$$
$$+ \alpha_n \langle (\gamma f - \mu G) x^*, j_q(x_{n+1} - x^*) \rangle.$$

This implies that

$$\|x_{n+1} - x^*\|^q \leq \frac{1 - \alpha_n(\tau - \gamma\beta)}{1 + \alpha_n(q-1)(\tau - \gamma\beta)} \|x_n - x^*\|^q$$
$$+ \frac{q\alpha_n}{1 + \alpha_n(q-1)(\tau - \gamma\beta)} \langle (\gamma f - \mu G) x^*, j_q(x_{n+1} - x^*) \rangle$$
$$\leq [1 - \alpha_n(\tau - \gamma\beta)] \|x_n - x^*\|^q$$
$$+ \frac{q\alpha_n}{1 + \alpha_n(q-1)(\tau - \gamma\beta)} \langle (\gamma f - \mu G) x^*, j_q(x_{n+1} - x^*) \rangle$$
$$\leq (1 - \gamma_n) \|x_n - x^*\|^q + \delta_n,$$

where $\gamma_n := \alpha_n(\tau - \gamma\beta)$ and $\delta_n := \frac{q\alpha_n}{1+\alpha_n(q-1)(\tau-\gamma\beta)} \langle (\gamma f - \mu G)x^*, j_q(x_{n+1} - x^*) \rangle$. From (K1′), $\lim_{n\to\infty} \gamma_n = 0$, $\sum_{n=0}^{\infty} \gamma_n = \infty$. $\frac{\delta_n}{\gamma_n} = \frac{q}{[1+\alpha_n(q-1)(\tau-\gamma\beta)](\tau-\gamma\beta)} \langle (\gamma f - \mu G)x^*, j_q(x_{n+1} - x^*) \rangle$. So, $\limsup_{n\to\infty} \frac{\delta_n}{\gamma_n} \leq 0$. Hence, by Lemma 2.8 we have that $x_n \to x^*$ as $n \to \infty$. This completes the proof. □

## 5 Conclusion

Let $E = H$ be a real Hilbert space, $q = 2$, $d_q = 1$ in Theorems 3.1 and 4.1, then we get the following result.

**Corollary 5.1** (Tian and Di [12]) *Let $\{x_n\}$ be a sequence generated by*

$$\begin{cases} T^{\beta_n} = \beta_n I + (1-\beta_n) \sum_{i=1}^{N} \lambda_i T_i, \\ x_{n+1} = \alpha_n \gamma f(x_n) + (I - \alpha_n \mu G) T^{\beta_n} x_n, \quad n \geq 0. \end{cases}$$

*Assume that $\{\alpha_n\}$ and $\{\beta_n\}$ are sequences in $(0,1)$ satisfying the conditions*

(K1) $\lim_{n\to\infty} \alpha_n = 0$, $\sum_{n=1}^{\infty} \alpha_n = \infty$,
(K2) $\sum_{n=1}^{\infty} |\alpha_{n+1} - \alpha_n| < \infty$, $\sum_{n=1}^{\infty} |\beta_{n+1} - \beta_n| < \infty$,
(K3) $0 < \max k_i \leq \beta_n < a < 1$, $\forall n \geq 0$,

*then $\{x_n\}$ converges strongly to a common fixed point of $\{T_i\}_{i=1}^{N}$, which solves the variational inequality*

$$\langle (\gamma f - \mu G)x^*, x - x^* \rangle \leq 0, \quad \forall x \in \bigcap_{i=1}^{N} F(T_i).$$

**Corollary 5.2** (Tian and Di [12]) *Let $\{x_n\}$ be a sequence generated by*

$$\begin{cases} A_{[n]} = \beta_{[n]} I + (1 - \beta_{[n]}) T_{[n]}, \\ x_{n+1} = \alpha_n \gamma f(x_n) + (I - \alpha_n \mu G) A_{[n+1]} x_n, \quad n \geq 0, \end{cases}$$

where $T_{[n]} = T_i$, with $i = n(\mod N)$, $1 \leq i \leq N$. Assume that $\{\alpha_n\}$ and $\{\beta_n\}$ are sequences in $(0,1)$ satisfying the conditions:

(K1′) $\lim_{n\to\infty} \alpha_n = 0$, $\sum_{n=1}^{\infty} \alpha_n = \infty$;
(K2′) $\sum_{n=1}^{\infty} |\alpha_{n+1} - \alpha_n| < \infty$ or $\lim_{n\to\infty} \frac{\alpha_n}{\alpha_{n+N}} = 1$;
(K3′) $\beta_{[n]} \in [k, 1)$, $\forall n \geq 0$, where $k = \max\{k_i : 1 \leq i \leq N\}$.

Then $\{x_n\}$ converges strongly to a common fixed point of $\{T_i\}_{i=1}^{N}$, which solves the variational inequality

$$\langle (\gamma f - \mu G) x^*, x - x^* \rangle \leq 0, \quad \forall x \in \bigcap_{i=1}^{N} F(T_i).$$

Let $E = H$ be a real Hilbert space; $q = 2$, $d_q = 1$, $n = 1$, $\beta_n = 0$, $G = A$, $\mu = 1$ and $T$ is a nonexpansive mapping in Theorems 3.1 and 4.1, then we get the following.

**Corollary 5.3** (Tian [11]) *Let $\{x_n\}$ be a sequence generated by $x_0 \in H$,*

$$x_{n+1} = \alpha_n \gamma f(x_n) + (I - \mu \alpha_n F) T x_n, \quad n \geq 0.$$

*Assume that $\{\alpha_n\}$ is a sequence in $(0,1)$ satisfying the conditions:*

(C1) $\lim_{n\to\infty} \alpha_n = 0$,
(C2) $\sum_{n=0}^{\infty} \alpha_n = \infty$,
(C3) *either* $\sum_{n=1}^{\infty} |\alpha_{n+1} - \alpha_n| < \infty$ *or* $\lim_{n\to\infty} \frac{\alpha_{n+1}}{\alpha_n} = 1$,

*then $\{x_n\}$ converges strongly to a common fixed point $\tilde{x}$ of $T$, which solves the variational inequality*

$$\langle (\gamma f - \mu F) \tilde{x}, x - \tilde{x} \rangle \leq 0, \quad \forall x \in F(T).$$

Let $E = H$ be a real Hilbert space; $q = 2$, $d_q = 1$, $n = 1$, $\beta_n = 0$, $G = A$, $\mu = 1$ and $T$ is a nonexpansive mapping in Theorems 3.1 and 4.1, then we get the following.

**Corollary 5.4** (Marino and Xu [10]) *Let $\{x_n\}$ be a sequence generated by $x_0 \in H$,*

$$x_{n+1} = \alpha_n \gamma f(x_n) + (I - \alpha_n A) T x_n, \quad n \geq 0.$$

*Assume that $\{\alpha_n\}$ is a sequence in $(0,1)$ satisfying* (C1)-(C3), *then $\{x_n\}$ converges strongly to a common fixed point $\tilde{x}$ of $T$, which solves the variational inequality*

$$\langle (\gamma f - A) \tilde{x}, x - \tilde{x} \rangle \leq 0, \quad \forall x \in F(T).$$

Let $E = H$ be a real Hilbert space; $q = 2$, $d_q = 1$, $n = 1$, $\beta_n = 0$, $G = F$, $\alpha_n = \lambda_n$, $\gamma = 0$ and $T$ is a nonexpansive mapping in Theorems 3.1 and 4.1, then we get the following.

**Corollary 5.5** (Yamada [13]) *Let $\{x_n\}$ be a sequence generated by $x_0 \in H$,*

$$x_{n+1} = T x_n - \mu \lambda_n F(T x_n), \quad n \geq 0.$$

*Assume that $\{\lambda_n\}$ is a sequence in $(0,1)$ satisfying the conditions:*

(i) $\lim_{n\to\infty} \lambda_n = 0$,

(ii) $\sum_{n=1}^{\infty} \lambda_n = \infty$,

(iii) *either* $\sum_{n=1}^{\infty} |\lambda_{n+1} - \lambda_n| < \infty$ *or* $\lim_{n\to\infty} \frac{\lambda_{n+1}}{\lambda_n} = 1$, *then* $\{x_n\}$ *converges strongly to a common fixed point of* $T$, *which solves the variational inequality*

$$\langle F\tilde{x}, x - \tilde{x} \rangle \geq 0, \quad \forall x \in F(T).$$

**Corollary 5.6** (Yamada [13]) *Let* $\{x_n\}$ *be a sequence generated by* $x_0 \in H$,

$$x_{n+1} = T^{\lambda_n} x_n = (I - \mu \lambda_n F) T_{[n]} x_n.$$

*Assume that* $\{\lambda_n\}$ *is a sequence in* $(0,1)$ *satisfying* (i)-(iii), *then* $\{x_n\}$ *converges strongly to a common fixed point of* $T$, *which solves the variational inequality*

$$\langle F\tilde{x}, x - \tilde{x} \rangle \geq 0, \quad \forall x \in F(T).$$

**Author details**
[1]Department of Mathematical Sciences, College of Remedial and Advanced Studies, P.M.B. 048, Kafin Hausa, Jigawa, Nigeria. [2]Department of Mathematical Sciences, Bauchi State University, Gadau, Nigeria. [3]Department of Mathematics, Gombe State University, Gombe, Nigeria.